\documentclass[twoside]{article}
\def\nkyear{200?}
\def\nkvol{??}
\def\nkno{?}

\def\firstpage{1}
\def\lastpage{20}
\def\correction{4}
\def\shorttitle{Iteration formulae} 
\def\shortauthor{{\it L. Wu\ \ and\ \ C. Zhu}} 
\usepackage b \baselineskip 14.5pt

\input{cam.doi}
\journalnumber{11401}
\DOImsnr{0001}
 \DOIyear{2007}

\usepackage{url}

\theoremstyle{thmstyle}
\newtheorem{theorem}{\indent Theorem}[section]
\newtheorem{corollary}{\indent Corollary}[section]
\newtheorem{lemma}{\indent Lemma}[section]
\newtheorem{proposition}{\indent Proposition}[section]
\newtheorem{definition}{\indent Definition}[section]

\theoremstyle{remark}

\newcommand{\exendproof}{\renewcommand{\qed}{\relax}\end{proof}}

\newsavebox{\SmallMathBox}

\DeclareRobustCommand*{\nicefrac}[2]{\ifmmode\mathnicefrac{#1}
{ #2}%
  \else\textnicefrac{#1}{#2}\fi}
\newcommand*{\textnicefrac}[2]{\check@mathfonts%
\mbox{\raisebox{.5ex}{\fontsize\sf@size\z@\selectfont#1}\kern-.
1em%
/\kern-.1em\raisebox{- .25ex}{\fontsize\sf@size\z@\selectfont#2} }}
\newcommand*{\mathnicefrac}[2]{%
  \mathchoice
    {\m@fr@c{\scriptstyle}{#1}{#2}}
    {\m@fr@c{\scriptstyle}{#1}{#2}}
    {\m@fr@c{\scriptscriptstyle}{#1}{#2}}
    {\m@fr@c{\scriptscriptstyle}{#1}{#2}}}

\def\lla{\langle}

\def\noi{\noindent}

\def\rra{\rangle}
\def\sqm1{\sqrt{-1}}

\def\wt{\widetilde}

\def\={\cong}
\def\>{\supset}
\def\<{\subset}

\def\12{\frac{1}{2}}
\def\0{^{\circ}}

\def\PPP{{\rm P}}
\def\SSS{{\rm S}}
\def\UUU{{\rm U}}

\def\CC{{\mathbb C}}

\def\RR{{\mathbb R}}

\def\ZZ{{\mathbb Z}}

\def\Bb{{\mathcal B}}

\def\Ff{{\mathcal F}}

\def\Kk{{\mathcal K}}
\def\Ll{{\mathcal L}}

\def\Pp{{\mathcal P}}

\def\a{\alpha}

\def\C{\CC}

\def\f{\varphi}

\def\K{\Kk}
\def\la{\lambda}

\def\n{\nu}

\def\R{\RR}

\def\w{\omega}

\def\Z{\ZZ}

\newcommand{\defeq}{\mathrel{\mathop:}=}

\DeclareMathOperator{\End}{End}

\DeclareMathOperator{\GL}{GL}  \DeclareMathOperator{\Graph}{Gr}
 
\DeclareMathOperator{\Hom}{Hom}

\DeclareMathOperator{\Index}{index}

\DeclareMathOperator{\Mas}{Mas}

 \DeclareMathOperator{\ran}{im}

 \DeclareMathOperator{\Sp}{Sp} \DeclareMathOperator{\ssp}{sp}

\title{\Large \bf \boldmath\ \\ Iteration formulae of the Maslov-type index theory in weak symplectic Hilbert spaces$^{\ast}$} 

\author{\large Li Wu$^1$\qquad Chaofeng Zhu$^{2}$} 

\date{}

\begin{document}

\maketitle

\thispagestyle{first}
\renewcommand{\thefootnote}{\fnsymbol{footnote}}

\footnotetext{\hspace*{-5mm} \begin{tabular}{@{}r@{}p{13.4cm}@{}}
& Manuscript received January 30, 2016\\ 
$^1$ & Chern Institute of Mathematics and LPMC,
Nankai University,
Tianjin 300071, P. R. China.\\
&{E-mail:nankai.wuli@gmail.com} \\
$^{2}$ & Chern Institute of Mathematics and LPMC,
Nankai University,
Tianjin 300071, P. R. China.\\
&{E-mail:zhucf@nankai.edu.cn}\\
&Corresponding author: CZ [\texttt{zhucf@nankai.edu.cn}]\\
 $^{\ast}$ & Project supported by
  the National Natural Science Foundation of China (No.11221091 and No. 11471169) and LPMC of MOE of China.
\end{tabular}}

\renewcommand{\thefootnote}{\arabic{footnote}}

\begin{abstract} 
 In this paper, we prove a splitting formula for the Maslov-type indices of symplectic paths induced by the splitting of the nullity in weak symplectic Hilbert space. Then we give a direct proof of iteration formulae for the Maslov-type indices of symplectic paths.

\vskip 4.5mm

\nd \begin{tabular}{@{}l@{ }p{10.1cm}} {\bf Keywords } &
Maslov-type index, positive path, iteration formula
\end{tabular}

\nd {\bf 2010 MR Subject Classification } 
53D12, 58J30

\end{abstract}

\baselineskip 14pt

\setlength{\parindent}{1.5em}

\setcounter{section}{0}

\section{Introduction}\label{s:introduction}

Since the pioneer work \cite{Bo56} on iteration of closed geodesics, the theory was generalized to convex Hamiltonian case \cite{Ek90}, to symplectic paths \cite{Lo99}, and to brake orbits \cite{LiZh14a}. The iteration theory for symplectic paths turns out to be extremely useful in the study of the periodical solutions of the finite dimensional Hamiltonian systems \cite{Ek90,SaZe92,Lo02,LoZh02,LoZhZh06,Hi09,Gi10,LiZh14a,LiZh14}.

All the known proofs of the iteration formulae follow the following line. Firstly, they give the relation between the (relative) Morse index for a certain differential operator and the Maslov-type index of its fundamental solution. Then they show that such an operator is in diagonal form under the iteration with respect to a certain orthogonal decomposition of the subspaces. The iteration formulae then follow from the above argument.

In this paper, we shall give a direct proof of the iteration formulae. The idea is that each splitting formula of the nullity will induce a splitting formula of the Maslov-type indices. The iteration formula then follows from the homotopy invariance of the Maslov-type indices.

Our method has three advantages. Firstly, we do not require that the symplectic path starts from the identity. So we can treat the iteration theory for symplectic paths in symplectic Hilbert spaces. It will be used in the study of the periodical solutions of the infinite dimensional Hamiltonian systems. Secondly, we do not need the orthogonality of the decomposition. So we can get more general formulae in some cases. It can be applied to non-contractible periodical solutions of Hamiltonian systems. Finally, our results apply to symplectic paths on weak symplectic Hilbert spaces.

\section{A Splitting formula of the Maslov-type index induced by the splitting of the nullity}\label{s:nullity-Maslov-index}

In this paper we denote by $\dim V$ the complex dimension of a complex vector space $V$. A pair of complex linear subspaces $(\lambda,\mu)$ of $V$ is called {\em Fredholm} if $\dim(\lambda\cap\mu)$ and $\dim V/(\lambda+\mu)$ are finite. The {\em index} of the pair $(\lambda,\mu)$ in $V$ is defined as
\begin{align}\label{e:fred-index}\Index(\lambda,\mu)=\dim(\lambda\cap\mu)-\dim V/(\lambda+\mu).
\end{align}

Let $(H,\lla\cdot,\cdot\rra)$ be a complex Hilbert space with a bounded injective map $J$ such that $J^*=-J$. Then we have the symplectic structure $\omega$ on $H$ defined by $\omega(x,y):=\lla Jx,y\rra$ for each $x,y\in H$. Since $J^*$ is injective, $\ran J$ is dense in $H$. The operator $J^{-1}$ on $H$ is a closed operator with dense domain $\ran J$. The {\em annihilator} of a linear subspace ${\la}$ of $H$ is defined by
\[
{\la}^{\w} \ :=\  \{y\in X ; \w(x,y)\ =\ 0\text{ for all $x\in {\la}$}\} .
\]
The subspace $\la$ is called a {\em Lagrangian subspace} of $H$ if $\la=\la^{\w}$. We denote by $\Ll(H)$ the set of the Lagrangian subspaces of $H$.
We define
\begin{align}\label{e:FLH}
\Ff\Ll(H)\ :&=\ \{(\lambda,\mu)\in(\Ll(H))^2;  (\lambda,\mu) \text{ is Fredholm}\},\\
\label{e:FLkH}\Ff\Ll_k(H)\ :&=\ \{(\lambda,\mu)\in\Ff\Ll(H); \Index(\la,\mu)=k\}.
\end{align}

Denote by $\GL(H)$, $\Bb(H)$, $\SSS(H)$, $\PPP(H)$ and $\UUU(H)$ the set of bounded invertible operators,
bounded operators, bounded selfadjoint operators, positive definite operators and unitary operators on $H$ respectively. The {\em symplectic group} $\Sp(H)$ and the {\em symplectic Lie algebra} $\ssp(H)$ are defined as
\begin{align}\label{e:symplectic-graoup-H}\Sp(H):&=\{M\in \GL(H);M^*JM=J\},\\
\label{e:symplectic-Lie-algebra-H}\ssp(H):&=\{M\in\Bb(H);M^*J+JM=0\}.
\end{align}
Then we have $J\in\ssp(H)$ and $e^{Js}\in\Sp(H)$ for $s\in\R$. The path $e^{Js}$, $s\in\R$ is crucial in our main result, Theorem \ref {t:nullity-Maslov-index} below. So our method does not work in the general symplect Banach space case. For each $M\in\Bb(H)$ such that $M^*JM=J$, the operator $M$ is injective. If $\dim H$ is infinite and $J$ is invertible, there exists an $M\in\Bb(H)\setminus\GL(H)$ such that $M^*JM=J$. For each $M\in\Sp(H)$, we have $M^{-1}=J^{-1}M^*J$ and $\ran (M^*J)\subset\ran J$. Let $z\in\C$ and $|z|=1$. Then we have $zI_H\in\Sp(H)$, where $I_H$ denotes the identity map on $H$. The fact is used in Definition \ref{d:Maslov-type-index} below. So we use the concept of complex symplectic Hilbert space insead of the real one.

Let $\lambda(s)$, $s\in(-\varepsilon,\varepsilon)$ be a path in $\Ll(H)$ differentiable at $s=0$ and $\lambda(0)$ be a {\em L-complemented Lagrangian subspace}, i.e., there exists a $\mu\in\Ll(H)$ such that $H=\lambda(0)\oplus\mu$. If $J$ is an isomorphism, each Lagrangian subspace of $H$ is a L-complemented Lagrangian subspace. Define the form $q(\lambda,0)$ on $\lambda(0)$ by
\begin{align}\label{e:form-q}
q(\lambda,0)(x,y):=\frac{d}{ds}|_{s=0}\omega(x,y(s)),
\end{align}
where $x,y\in\lambda(0)$, $y(s)\in\lambda(s)$, and $y(s)-y\in\mu$. Then the form $q(\lambda,0)$ is independent on the choices of $\mu\in\Ll(H)$ with $H=\lambda(0)\oplus\mu$. The path $\lambda$ is said to be {\em (semi-)positive} (cf. \cite[Definition 3.2.8]{BoZh14}) if $q(\lambda,s)$ is (semi-)positive definite for each $s\in(-\varepsilon,\varepsilon)$.

\begin{lemma}\label{l:gamma-crossing-form} Let $\gamma(s)$, $s\in(-\varepsilon,\varepsilon)$ be a path in $\Sp(H)$ differentiable at $s=0$ with $\gamma(0)=I_H$ and $\lambda(0)$ be a L-complemented Lagrangian subspace. Set $\lambda(s):=\gamma(s)\lambda(0)$ for each $s\in(-\varepsilon,\varepsilon)$. Then we have
\begin{equation}\label{e:gamma-crossing-form}
q(\lambda,0)(x,y)=\frac{d}{ds}|_{s=0}\omega(x,\gamma(s)y),
\end{equation}
where $x,y\in\lambda(0)$.
\end{lemma}

\begin{proof} Let $\mu\in\Ll(H)$ be such that $H=\lambda(0)\oplus\mu$. Let $\gamma(s)$ be of the form
\[\gamma(s)=\left(\begin{array}{cc}A(s) & B(s)\\C(s) & D(s) \end{array}\right).\]
For $x,y\in\lambda(0)$, let $y(s)\in\lambda(s)$ be such that
$y(s)-y\in\mu$. Then we have $y(s)=y+C(s)A(s)^{-1}y$, and $\omega(x,y(s)-\gamma(s)y)=\omega(x,C(s)(A(s)^{-1}-I_{\lambda(0)})y)$. By (\ref{e:form-q}),
our result (\ref{e:gamma-crossing-form}) holds.
\end{proof}

Set $X:=H\times H$.
We define the symplectic structure of $X$ by
\begin{equation}
 \tilde \omega(v,w)=\langle\tilde J v,w\rangle,\forall v,w\in X, \text{ where } \tilde J=(-J)\oplus J.
\end{equation}
For any $M\in \Sp(H)$, we have $X=\Graph(M)\oplus\Graph(-M)$, where
\[\Graph(M)\defeq \{(x,Mx),x\in H\}.\]
By \cite[Lemma 4]{BoZh13}, $\Graph(M)$ is a L-complemented Lagrangian subspace of $X$.

\begin{definition}\label{d:Maslov-type-index}(cf. \cite[Definition 4.6]{Zh06}) Let $(H,\omega)$ and $(X,\tilde \omega)$ be symplectic Hilbert spaces as defined above. Let $V$ be a Lagrangian subspace of $X$. Let $\gamma(t)$, $t\in[a,b]$ be a path in $\Sp(H)$. Assume that $(\Graph(\gamma(t)),V)$ is a Fredholm pair for each $t\in[a,b]$. Denote by $\Mas\{\cdot,\cdot\}$ the Maslov index for the path of Fredholm pairs of Lagrangian subspaces in $X$ with index $0$ defined by \cite[Definition 3.1.4]{BoZh14}. The {\em Maslov-type index} $i_V(\gamma)$ is defined by
\begin{equation}\label{e:Maslov-type-index}i_V(\gamma)=\Mas\{\Graph(\gamma),V\}.\end{equation}
For each $M\in\Sp(H)$, we define the {\em nullity} $\nu_V(M)$ and {\em co-nullity} $\tilde\nu_V(M)$ by
\begin{equation}\label{e:nullity}\nu_V(M)=\dim(\Graph(M)\cap V),\quad\tilde\nu_V(M)=\dim(X/(\Graph(M)+V).\end{equation}
If $N\in\Sp(H)$ and $V=\Graph(N)$, we define $i_N(\gamma):=i_V(\gamma)$, $\nu_N(M):=\nu_V(M)$ and $\tilde\nu_N(M):=\tilde\nu_V(M)$. If $z\in\C$, $|z|=1$ and $V=\Graph(z I_H)$, we define $i_z(\gamma):=i_V(\gamma)$, $\nu_z(M):=\nu_V(M)$ and $\tilde\nu_z(M):=\tilde\nu_V(M)$.
\end{definition}

By \cite[Lemma 4.4]{Zh06}, for each $N\in\Sp(H)$ and $\gamma\in C([a,b],\Sp(H))$ we have
\begin{equation}\label{e:N-Maslov-type}i_N(\gamma)=i_1(\gamma N^{-1})=i_1(N^{-1}\gamma).\end{equation}
By \cite[Proposition 1]{BoZh13}, for each $V\in\Ll(X)$ and $M\in\Sp(H)$, we have
\begin{equation}\label{e:co-nullity}\nu_V(M)\le\tilde\nu_V(M).
\end{equation}

\begin{definition}\label{d:positive-path}(cf. \cite[Definition 3.1]{Bo56}) Let $(H,\omega)$ be a symplectic Hilbert space as defined above.
A $C^1$ path $\gamma:[a,b]\rightarrow \Sp(H)$ is called a {\em positive path} if $\dot\gamma(t)\gamma(t)^{-1}\in K$ for each $t\in [a,b]$, where we denote by $\dot\gamma:=\frac{d \gamma}{dt}$, and the cone $K$ is defined by
\[K=\{A\in\Bb(H);-JA=A^*J>0\}.\]
\end{definition}

By \cite[Lemma 3.1]{Du76}, the path $\gamma$ is positive if and only if $\Graph(\gamma)$ is positive.

The following lemma is our key observation.

\begin{lemma}\label{l:product-positive-path} Let $(H,\omega)$ be a symplectic Hilbert space as defined above.
Let $\gamma_1(t)$, $\gamma_2(t)$ with $t\in[a,b]$ be two positive paths in $\Sp(H)$.
Then $\gamma(t):=\gamma_1(t)\gamma_2(t),t\in [a,b]$ is a positive path.
\end{lemma}

\begin{proof}
We have
\begin{align*}
-J\dot\gamma\gamma^{-1}&=-J(\dot\gamma_1\gamma_2+\gamma_1\dot\gamma_2)(\gamma_1\gamma_2)^{-1}\\
&=-J\dot\gamma_1\gamma_1^{-1}+(\gamma_1^{*})^{-1}(-J\dot\gamma_2\gamma_2^{-1})(\gamma_1)^{-1}.
\end{align*}
Since $\gamma_1$, $\gamma_2$ are both positive paths, $-J\dot\gamma_1\gamma_1^{-1}$ and
$(\gamma_1^{*})^{-1}(-J\dot\gamma_2\gamma_2^{-1})(\gamma_1)^{-1}$ are both positive definite.
So $-J\dot\gamma\gamma^{-1}$ is positive definite and the lemma follows.
\end{proof}

We need a special case of \cite[Theorem 3.2.12]{BoZh14}.

\begin{theorem}\label{t:comparision}
Let $(H,\omega)$ and $ (\wt{H},\wt{\omega})$ be two symplectic Hilbert space.
For $0\le a\le \delta$, $\delta>0$, we are given continuous two-parameter families
\begin{align}\label{e:two-parameter}
&\{(\la(s,a),\mu(s))\in (\Ll(H,\w))^2\}\text{ and}\\
\nonumber&\{(\wt\la(s,a),\wt\mu(s))\in(\Ll(\wt H,\wt\w))^2\}.
\end{align}
We assume that
\begin{align}
\label{e:comparision0} &(\la(s,0),\mu(s))\in \Ff\Ll_0(H)\text{ and }
(\wt\la(s,0),\wt\mu(s))\in\Ff\Ll_0(\wt H),\\
\label{e:comparisionI} &\{\la(s,a)\} \text{ differentiable in $a$ and semi-positive for fixed $s$},\\
\label{e:comparisionII} &\{\wt\la(s,a)\} \text{ differentiable in $a$ and positive for fixed $s$},\\
\label{e:comparisionIII} &\dim(\la(s,a)\cap\mu(s)) - \dim(\wt\la(s,a)\cap\wt\mu(s)) = c(s).
\end{align}
Then we have
\begin{equation}\label{e:comparisionIV}
\Mas\{\la(s,0),\mu(s);\w\} = \Mas\{\wt\la(s,0),\wt\mu(s);\wt\w\} .
\end{equation}
\end{theorem}

The main result of the paper is following. Note that the path $Me^{Js}$, $s\in\R$ is positive for each $M\in\Sp(H)$.

\begin{theorem}[Splitting of the Maslov-type index induced by splitting of the nullity]\label{t:nullity-Maslov-index}
Let $(H,\omega)$ and $(X,\tilde \omega)$ be symplectic Hilbert spaces as defined above.
Let $f_j:\Sp(H)\rightarrow \Sp(H)$ ($j=1,\ldots,k$) be a family of $C^1$ maps, where $k$ is a positive integer.
Assume that there are Lagrangian subspaces $\{V_j\}_{j=0,\ldots,k}$ of $X$ such that the following hold for each $M\in \Sp(H)$.
\begin{itemize}
\item[(i)] the pair $(\Graph(\prod_{1\leq j\leq k} f_k(M)),V_0)$ is a Fredholm pair with index $0$ if and only if all pairs $(\Graph (f_j(M)),V_j)$ are Fredhlom pairs with index $0$, and
\begin{equation}
 \nu_{V_0}(\prod_{1\leq j\leq k} f_j(M))=\sum_{1\leq j\leq k}  \nu_{V_j} (f_j(M)).
\end{equation}
\item[(ii)]
$f_j( e^{Js}M )$, $s\in\R$ is positive for each $j=1,\ldots, k$.
\end{itemize}
Let $\gamma\in C([a,b],H)$ such that the pair $(\Graph(\prod_{1\leq j\leq k} f_j(\gamma(t))),V_0)$ is a Fredholm pair with index $0$ for each $t\in[0,1]$.
Then we have
\begin{equation}\label{e:nullity-Maslov-index}
 i_{V_0}(\prod_{1\leq j\leq k} f_j(\gamma))=\sum_{1\leq j\leq k}i_{V_j} (f_j(\gamma)).
\end{equation}
\end{theorem}

\begin{proof} Let $X:=X$ and $\tilde X:=X^k$. For each $(s,t)\in\R\times[0,1]$, set
\begin{align*}
\la(t,s):=\Graph(\prod_{1\leq j\leq k} f_j(e^{Js}\gamma(t))),\quad &\mu(t):=V_0,\\
\tilde\la(t,s):=\prod_{1\leq j\leq k}\Graph(f_j(e^{Js}\gamma(t))),\quad &\tilde\mu(t):=\prod_{1\leq j\leq k}V_j.
\end{align*}
By (i), for each $t\in[0,1]$, there hold that the pair $(\la(t,0),\mu(t))$ is a Fredholm pair of Lagrangian subspace in $X$,  $(\tilde\la(t,0),\tilde\mu(t))$ is a Fredholm pair of Lagrangian subspace in $\tilde X$, and
$\dim(\la(t,s)\cap\mu(t))=\dim(\tilde\la(t,s)\cap\tilde\mu(t))$. By (ii) and Lemma \ref{l:product-positive-path}, the two paths $\la(t,\cdot)$ and $\tilde\la(t,\cdot)$ are both positive for each fixed $t$.
By \cite[Theorem 2.2.1.c]{BoZh14}, the Maslov index is additive under direct sum.
We have
\[
 \Mas\{\tilde\la(\cdot,0),\tilde\mu;\tilde X\}=\sum_{1\le j\le k}\Mas\{\Graph(f_j(\cdot,0)),V_j\}
\]

By Theorem \ref{t:comparision} we have
\begin{align*}
i_{V_0}(\prod_{1\leq j\leq k} f_j(\gamma))&=\Mas\{\la(\cdot,0),\mu;X\}=\Mas\{\tilde\la(\cdot,0),\tilde\mu;\tilde X\}\\
&=\sum_{1\le j\le k}\Mas\{\Graph(f_j(\cdot,0)),V_j\}=\sum_{1\leq j\leq k}i_{V_j} (f_j(\gamma)).
\end{align*}
\end{proof}

\begin{rem}
Each pair of Lagrangian subspaces of a finite dimensional symplectic space is a Fredhlom pair with index $0$.
We will choose suitable functions $\{f_i\}$ to get some classical iteration formulae for the symplectic paths in section \ref{s:iteration}.
\end{rem}

\section{The global structure of $\Sp(H)$ when $\omega$ is strong}\label{s:globle}

Let $(H,\omega)$ be the symplectic Hilbert space defined in \S\ref{s:nullity-Maslov-index}. Assume that $\omega$ is strong, i.e., $J$ is an isomorphism. In this section we give the global structure of $\Sp(H)$.

\begin{lemma}\label{l:normalize-J} Let $(H,\w,J)$ be a strong symplectic Hilbert space. Set $J_1:=(-J^2)^{-\frac{1}{2}}J$ and $\w_1(x,y):=\lla J_1x,y\rra$ for each $x,y\in H$. Then the following hold.
\newline (a) $(H,\w_1,J_1)$ is a strong symplectic Hilbert space, $J_1^2=-I_H$, $J_1^*=-J_1$.
\newline (b) $(-J^2)^{\frac{1}{4}}:(H,\omega)\to (H,\omega_1)$ is symplectic.
\newline (c) We have a homeomorphism $\varphi:\Sp(H,\w)\to\Sp(H,\w_1)$ defined by $\varphi(M)=(-J^2)^{\frac{1}{4}}M(-J^2)^{\frac{1}{4}}$ for each $M\in\Sp(H,\w)$.
\end{lemma}

\begin{proof} Direct computation. Note that $J(-J^2)^t=(-J^2)^tJ$ for each $t\in\R$.\end{proof}

By the lemma above, we can assume that $J^2=-I_H$. By Kuiper's theorem, $\UUU(H)$ is contractible if $\dim H=+\infty$. By the method similar to that in \cite[\S1.1,\S2.2]{Lo02}, we get the following Proposition \ref{p:structure-Sp}. Here we give a short proof.

\begin{proposition}\label{p:structure-Sp} Let $(H,\omega)$ be a symplectic Hilbert space defined in \S\ref{s:nullity-Maslov-index}. Assume that $J^2=-I_H$. Set $H^{\pm}:=\dim\ker(J\mp\sqrt{-1})$. Then $H=H^+\oplus H^-$.
\newline (a) Each $M\in\Sp(H)$ is uniquely represented by $M=AU$, where $A\in\PPP(H)\cap\Sp(H)$ and $U\in\UUU(H)\cap\Sp(H)$. So we have $\Sp(H)\cong (\PPP(H)\cap\Sp(H))\times (\UUU(H)\cap\Sp(H))$.
\newline (b) Each $A\in\PPP(H)\cap\Sp(H)$ is uniquely represented by $A=\exp(S)$, where $S\in\SSS(H)\cap\ssp(H)$ is of the form $S=\left(\begin{array}{cc}0 & S_{12} \\S_{12}^*& 0 \\ \end{array}\right)$, $S_{12}\in\Bb(H^-,H^+)$. So we have $\PPP(H)\cap\Sp(H)\cong \Bb(H^-,H^+)$.
\newline (c) Each $U\in\UUU(H)\cap\Sp(H)$ is of the form $U=\left(\begin{array}{cc}U_{11} & 0 \\0 & U_{22} \\ \end{array}\right)$, where $U_{11}\in\UUU(H^+)$, $U_{22}\in\UUU(H^-)$. So we have $\UUU(H)\cap\Sp(H)\cong \UUU(H^+)\times\UUU(H^-)$.
\newline (d) $\Sp(H)\cong \Bb(H^-,H^+)\times\UUU(H^+)\times\UUU(H^-)$ is path connected.
\newline (e) $\Sp(H)$ is path-connected, and the fundamental group of $\Sp(H)$ is given by
\begin{equation}\label{e:pi-1-Sp}\pi_1(\Sp(H))=\begin{cases}
0,&\text{if }\dim H^+=\dim H^-=0\text{ or }+\infty,\\
\Z,&\text{if }0<\dim H^{\pm}<+\infty,\dim H^{\mp}=0\text{ or }+\infty,\\
\Z\oplus\Z,&\text{if }0<\dim H^+,\dim H^-<+\infty.
\end{cases}
\end{equation}
\newline (f) $\Sp(H)$ is a real analytic Hilbert manifold modelled on the real Hilbert space $\ssp(H)$. Its real dimension is $m^2$ if $H=\C^m$.
\end{proposition}

\begin{proof} (a) By the polar decomposition, each $M\in\Sp(H)$ is uniquely represented by $M=AU$, where $A\in\PPP(H)$ and $U\in\UUU(H)$. Then we obtain
\[M=J^{-1}(M^*)^{-1}J=J^{-1}A^{-1}JJ^{-1}(U^*)^{-1}J.\]
By the uniqueness of the decomposition we obtain $A=J^{-1}AJ$ and $U=J^{-1}(U^*)^{-1}J$. Thus both $A$ and $U$ are symplectic. The converse is obvious.
\newline (b) Each $A\in\PPP(H)\cap\Sp(H)$ is uniquely represented by $A=\exp(S)$, where $S\in\SSS(H)$. Since $A$ is symplectic,
\[\exp(S)=J^{-1}(\exp(S))^{-1}J=\exp(-J^{-1}SJ).\]
By the uniqueness we obtain $S=-J^{-1}SJ$ and $S\in\ssp(H)$. The converse is obvious.
\newline (c) Note that $U\in\UUU(H)\cap\Sp(H)$ if and only if $UJ=JU$ and $U\in\UUU(H)$. Our results then follow from direct calculation.
\newline (d) follows from (a), (b), (c), and (e) follows from (d).
\newline (f) The first part of (f) follows from (d). If $H=\C^m$, we denote by $m^{\pm}:=\dim H^{\pm}$. Then the real dimension of $\Sp(H)$ is
\[2m^+m^-+(m^+)^2+(m^-)^2=(m^++m^-)^2=m^2.\]
\end{proof}

\section{Iteration formulae for symplectic paths}\label{s:iteration}

In this section we derive some iteration formulae for symplectic paths from Theorem \ref{t:nullity-Maslov-index}.

\subsection{Some general facts on Maslov-type indices in the finite dimensional case}\label{ss:facts}

The following lemmas is useful for deriving iteration formulae from Theorem \ref{t:nullity-Maslov-index}. For each $\tau>0$, we define
\begin{equation}\label{e:path-I-space}
\Pp_{\tau}(H):=\{\gamma\in C([0,\tau],\Sp(H));\gamma(0)=I_H\}.
\end{equation}

\begin{lemma}\label{l:Maslov-path-endpoints} Let $(H,\w)$ and $(X,\tilde\w)$ be symplectic Hilbert spaces defined in \S\ref{s:nullity-Maslov-index}. Assume that $\dim H<+\infty$. Let $\phi:\Pp_{\tau_1}(H)\to\Pp_{\tau_2}(H)$ be a continuous map such that $\phi(c_0)=c_0$, where $\tau_1$, $\tau_2$ are two positive numbers and $c_0$ denotes the constant path $I_H$. Let $V$ be a Lagrangian subspace of $X$.
Denote the curve $\gamma(s\cdot):[0,1]\rightarrow \Sp(H)$ by $\gamma(s\cdot)(t):=\gamma(st)$.
Then for each $\gamma\in\Pp_{\tau_1}(H)$, we have
\begin{equation}\label{e:Maslov-path-endpoints}i_V(\phi(\gamma))=i_V(\{\phi(\gamma(s\cdot))(\tau_2);s\in[0,1]\}).
\end{equation}
\end{lemma}

\begin{proof} Define the homotopy $h(s,t):[0,1]\times[0,\tau_2]\to\Sp(H)$ by $h(s,t):=\phi(\gamma(s\cdot))(t)$. Then we have $h(0,t)=h(s,0)=I_H$. By the homotopy invariance of the Maslov index, the equation (\ref{e:Maslov-path-endpoints}) follows.
\end{proof}

\begin{lemma}\label{l:Maslov-depend-endpoints} Let $(H,\w)$ and $(X,\tilde\w)$ be symplectic Hilbert spaces defined in \S\ref{s:nullity-Maslov-index}. Let $f_1,f_2:\Sp(H)\to\Sp(H)$ be a two continuous maps such that $f_1(I_H)=f_2(I_H)=I_H$. Let $V_1$, $V_2$ be two Lagrangian subspaces of $X$. Assume that $\dim H<+\infty$. Then the following hold.
\newline (a) Assume that $(f_1)_*=(f_2)_*:\pi_1(\Sp(H))\to\Z$. Then there is an interger valued function $\delta:\Sp(H)\to\Z$ defined by
\begin{equation}
\label{e:delta-fuction1}\delta(M):=i_{V_2}(f_2(\gamma))-i_{V_1}(f_1(\gamma)),
\end{equation}
where $\gamma\in\Pp_{\tau}(H)$, $\tau>0$ with $\gamma(\tau)=M$.
\newline (b) Assume that $(f_1)_*=(f_2)_*:\pi_1(\Sp(H))\to\Z$. Let $\delta:\Sp(H)\to\Z$ be the funtion defined in (a). Then for each $\gamma\in C([a,b],\Sp(H))$, we have
\begin{equation}
\label{e:delta-fuction2}\delta(\gamma(b))-\delta(\gamma(a))=i_{V_2}(f_2(\gamma))-i_{V_1}(f_1(\gamma)).
\end{equation}
\newline (c) We have $(f_1f_2)_*=(f_1)_*+(f_2)_*:\pi_1(\Sp(H))\to\Z$.
\end{lemma}

\begin{proof} (a) Let $\gamma_1$ and $\gamma_2$ be in $\Pp_{\tau}(H)$ with $\gamma_1(\tau)=\gamma_2(\tau)=M$. Let $\gamma_3\in \Pp_{2\tau}(H)$ be the loop in $\Sp(H)$ defined by $\gamma_1$ followed by the reverse of $\gamma_2$. Then we have $[f_1(\gamma_3)]=[f_2(\gamma_3)]\in\pi_1(\Sp(H))$. Since $\Ll(X)$ is path connected, there is a path  $\lambda:[0,1]\to\Ll(X)$ of Lagrangian subspaces such that $\lambda(0)=V_1$, $\lambda(1)=V_2$. Then we have a homotopy $h:[0,1]\times[0,2\tau]\to(\Ll(X))^2$ defined by $h(s,t)=(\lambda(s),\Graph(\gamma_3(t)))$. By the homotopy invariance of the Maslov index, we have $i_{V_1}(f_1(\gamma_3))=i_{V_2}(f_1(\gamma_3))=i_{V_2}(f_2(\gamma_3))$. Since $i_{V_j}(f_j(\gamma_3))=i_{V_j}(f_j(\gamma_1))-i_{V_j}(f_j(\gamma_2))$ for $j=1,2$, the equation (\ref{e:delta-fuction1}) follows.
\newline (b) Let $\gamma$ be in $C([a,b],\Sp(H))$. By Proposition \ref{p:structure-Sp}.d, there is a path $\alpha\in\Pp_{\tau}(H)$ such that $\alpha(\tau)=\gamma(a)$. Denote by $\tilde\gamma\in\Pp_{2\tau}(H)$ the path defined by $\alpha$ followed by $\gamma$. Since $i_{V_j}(f_j(\tilde\gamma))=i_{V_j}(f_j(\alpha))+i_{V_j}(f_j(\gamma))$ for $j=1,2$, the equation (\ref{e:delta-fuction2}) follows.
\newline (c) Let $\gamma\in C([0,1],\Sp(H))$ be such that $\gamma(0)=\gamma(1)=I$. Then we have a homotopy $h:[0,1]^2\to\Sp(H)$ such that $h(s,t)=f_1(\gamma)f_2(\gamma)$. So we have $[f_1(\gamma)f_2(\gamma)]=[f_1(\gamma)]+[f_2(\gamma)]$ in $\pi_1(\Sp(H))$ and (c) follows.
\end{proof}

\subsection{The iteration formulae for the $A$-iteration of the symplectic path}\label{ss:A-iteration}


Firstly we define the iteration of a symplectic path $\gamma\in\Pp_{\tau}(H)$.

\begin{definition}\label{d:A-iteration} (cf. \cite[(4.3),(4.4)]{LiZh14}) Let $(H,\omega)$ be a symplectic Hilbert space defined in \S\ref{s:nullity-Maslov-index}. Given an $A\in\Sp(2n)$, a $\tau>0$, a positive integer $k$, and a path $\gamma\in \Pp_{\tau}(H)$, we define \begin{equation}\label{e:A-iteration1}
\tilde\gamma(t)=A^j \gamma(t-j\tau)(A^{-1}\gamma(\tau))^j,\quad t\in [j\tau,(j+1)\tau],\;j=0,\ldots,k-1.
\end{equation}
\newline (a) We call $\tilde\gamma$ the {\em $k$-th $A$-iteration} of $\gamma$. The map $P_A(\gamma):=A^{-1}\gamma(\tau)$ is called the {\em Poincar\'{e} map} of $\gamma$ at $\tau$.
\newline (b) If $\dim H$ is finite, we define
\begin{equation}i_z(\gamma,k,A):=i_{zA^k}(\tilde\gamma),\quad \nu_z(\gamma,k,A):=\nu_{zA^k}(\tilde\gamma(k\tau))
\end{equation}
for each $z\in\C$ with $|z|=1$.
\end{definition}

We need the following lemmas.

\begin{lemma}\label{l:commutative-Fredholm} Let $X$ be a vector space and $A,B\in\End(X)$ be two linear maps. Assume that $AB=BA$. Then $AB$ is Freholm if and only $A$ and $B$ are Fredholm. In this case we have
\begin{equation}\label{e:product-Fredhom-index}\Index(AB)=\Index A+\Index B.
\end{equation}
\end{lemma}

\begin{proof} If $A$ and $B$ are Fredolm, $AB$ is Fredholm and (\ref{e:product-Fredhom-index}) holds. Since $AB=BA$, we have
\[\ker(AB)\supset \ker A+\ker B,\quad\ran(AB)\subset\ran A\cap\ran B.\]
So $A$ and $B$ are Fredholm if $AB$ is Fredholm.
\end{proof}

The following fact is clear.

\begin{lemma}\label{l:graph-Fredholm} Let $X$ and $Y$ be two abelian groups and $A,B\in\Hom(X,Y)$ be two group homomorphisms. Then we have $\ker(A-B)\cong\Graph(A)\cap\Graph(B)$ and $Y/\ran(A-B)\cong (X\times Y)/(\Graph(A)+\Graph(B))$.
\end{lemma}

We have the following form of iteration formula.

\begin{theorem}[The iteration formula for the power of the symplectic path]\label{t:iteration-power} Let $(H,\omega)$ be a symplectic Hilbert space defined in \S\ref{s:nullity-Maslov-index}. Let $\gamma\in C([a,b],\Sp(H))$ be a symplectic path such that $(\gamma(t))^k-I_H$ is Fredholm of index $0$ for each $t\in[a,b]$, where $k$ is a positive integer. Then we have
\begin{align}
\label{e:iteration-power}i_1( \gamma^k)&=\sum_{z^k=1}i_z(\gamma),\\
\label{e:power-nullity}\nu_1((\gamma(t))^k)&=\sum_{z^k=1}\nu_z(\gamma(t)),\quad \forall t\in[a,b].
\end{align}
\end{theorem}

\begin{proof} Set $V_j:=\Graph(\exp(\frac{2\pi j\sqrt{-1}}{k})I_H)$ for $j=0,\ldots k$. Recall that the integers $i_z(\gamma)$ and $\nu_z(\gamma(t))$ are defined by Definition \ref{d:Maslov-type-index}. It is clear that the equation (\ref{e:power-nullity}) holds.
By Lemma \ref{l:commutative-Fredholm} and \ref{l:graph-Fredholm}, we can apply Theorem \ref{t:nullity-Maslov-index} and our result follows.
\end{proof}

Let $\gamma(t)$ be the path in Theorem \ref{t:iteration-power}. Since the zero operator $0$ on $H$ is not Fredholm if $\dim H=+\infty$, $\gamma(t)$ can not be $I_H$ for each $t\in[a,b]$ in this situation.

Our next result for periodic Hermitian systems was obtained by \cite[Corollary of Theorem I]{Bo56} when $A=zI_H$ and $|z|=1$, in the real case was obtained by Y. Long \cite[Theorem 1.4]{Lo99} when $A=zI_H$ and $|z|=1$, by X. Hu and S. Sun \cite[Theorem 1.6]{HuSu09} when $A\in\Sp(2n,\R)\cap \UUU(2n)$, and by C. Liu and S. Tang \cite[Theorem 1.1]{LiTa15b} when $A\in\Sp(2n,\R)$.

\begin{corollary}[The iteration formula for the $A$-iteration of the symplectic path]\label{c:iteration-A} Let $(H,\omega)$ be a finite dimensional symplectic Hilbert space. Given an $A\in\Sp(H)$, a $\tau>0$, a positive integer $k$, and a path $\gamma\in \Pp_{\tau}(H)$, we have
\begin{align}
\label{e:index-iteration-A}i_1(\gamma,k,A)&=\sum_{z^k=1}i_z(\gamma,1,A),\\
\label{e:nullity-iteration-A}\nu_1(\gamma,k,A)&=\sum_{z^k=1}\nu_z(\gamma,1,A).
\end{align}
Moreover, there is a function $\delta_k:\Sp(H)\to\Z$ defined by
\[\delta_k(\gamma(\tau)):=i_1(\gamma,k,I_H)-ki_1(\gamma,1,I_H)\]
such that
\begin{equation}\label{e:iteration-function-A}
i_1(\gamma,k,A)-ki_1(\gamma,1,A)=\delta_k(A^{-1}\gamma(\tau))-\delta_k(A^{-1}).
\end{equation}

\begin{proof} Note that $A^{-k}\gamma(k\tau)=(A^{-1}\gamma(\tau))^k$. By (\ref{e:N-Maslov-type}), Lemma \ref{l:Maslov-path-endpoints} and Theorem \ref{t:iteration-power} we obtain (\ref{e:index-iteration-A}) and (\ref{e:nullity-iteration-A}). By Lemma \ref{l:Maslov-depend-endpoints}, the function $\delta_k$ is well-defined and (\ref{e:iteration-function-A}) holds.
\end{proof}
\end{corollary}

\subsection{Two times iteration formula for the generalized brake symmetry}\label{ss:2-iteration}

Firstly we define the iteration of a symplectic path $\gamma\in\Pp_{\tau}(H)$ for the brake symmetry.

\begin{definition}\label{d:brake-iteration} (cf. \cite[(4.3),(4.4)]{LiZh14}) Let $(H,\omega)$ be a symplectic Hilbert space defined in \S\ref{s:nullity-Maslov-index}. Let $N\in\Bb(H)$ be such that $N^*JN=-J$ and $N^2=I_H$. Given a $\tau>0$, a positive integer $k$, and a path $\gamma\in \Pp_{\tau}(H)$, we define \begin{equation}\label{e:brake-iteration1}
\gamma^{(k)}(t)=\begin{cases}\gamma(t-2j\tau)(\gamma(2\tau))^j,&t\in [2j\tau,(2j+1)\tau],j\in[0,\frac{k-1}{2}],j\in\Z,\\
N\gamma(2j\tau-t)N(\gamma(2\tau))^j,&t\in [(2j-1)\tau,2j\tau],j\in[1,\frac{k}{2}],j\in\Z,
\end{cases}
\end{equation}
where $\gamma(2\tau)=N\gamma(\tau)^{-1}N\gamma(\tau)$. We call $\gamma^{(k)}$ the {\em $k$-th $N$-brake iteration} of $\gamma$. The map $\gamma(2\tau)$ is called the {\em Poincar\'{e} map} of $\gamma$ at $2\tau$.
\end{definition}

We have the following algebraic results.

\begin{lemma}\label{l:index-pm}Let $X$, $Y$ be two linear spaces with decompositions $X=V_1\oplus V_2$ and $Y=W_1\oplus W_2$ respectively. Let $M\in\Hom(X,Y)$ be of the form $M=\left(\begin{array}{cc}A&B\\C&D\end{array}\right)$.
Then we have
\begin{align}
\label{e:ker-A}\ker A&=V_1\cap (M^{-1}W_2),\\
\label{e:Gr-MA}\Graph(M)+V_1\times W_2&=\left(\begin{array}{cc}I_{V_2}&0\\B&I_{W_1}\end{array}\right)(V_2\times\ran A)\oplus V_1\times W_2,\\
\label{e:index-pm}\Index A&=\Index(\Graph(M),V_1\times W_2).
\end{align}
\end{lemma}

\begin{proof} (\ref{e:ker-A}) is clear, (\ref{e:Gr-MA}) follows from direct calculation, and (\ref{e:index-pm}) follows from (\ref{e:ker-A}) and (\ref{e:Gr-MA}).
\end{proof}

\begin{lemma}\label{l:2S-iteration} Let $X$ be an abelian group such that $2I_X$ is an isomorphism on $X$. Let $M,N,S\in\End(X)$ be three endmorhpisms on $X$ such that $M$ is an isomorphsim and $N^2=(NS)^2=I_X$.
Set $U^{\pm}:=\ker(N\mp I_X)$, $V^{\pm}:=\ker(NS\mp I_X)$, and $K:=NM^{-1}NM$. Let $M$ be in the form
\begin{equation}\label{e:M-block-form}
M=\left(\begin{array}{cc}A&B\\C&D\end{array}\right):V^+\oplus V^-\to U^+\oplus U^-.
\end{equation}
Then we have
\begin{align}
\label{e:ker-ran-KS}\ker(K-S)=\ker C\oplus\ker B,&\quad\ran(MN(K-S))=\ran B\oplus\ran C,\\
\label{e:ker-C}\ker C=V^+\cap(M^{-1}U^+)&=\ker(K-S)\cap V^+,\\
\label{e:ker-B}\ker B=V^-\cap(M^{-1}U^-)&=\ker(K-S)\cap V^-.
\end{align}
\end{lemma}

\begin{proof} Since $(NS)^2=I_X$, we have
\begin{align*}(NM-MNS)(I_X\pm NS)&=NM(I_X\pm NS)-M(NS\pm I_X)\\
&=(NM\mp M)(I_X\pm NS)\\
&=(N\mp I_X)M(I_X\pm NS).
\end{align*}
Since $N^2=(NS)^2=I_X=\frac{I_X+N}{2}+\frac{I_X-N}{2}=\frac{I_X+NS}{2}+\frac{I_X-NS}{2}$, we have
\begin{equation}\label{e:NMS}
MN(K-S)=NM-MNS=\left(\begin{array}{cc}0&2B\\-2C&0\end{array}\right).
\end{equation}
Our results then follows from (\ref{e:NMS}) and direct computations.
\end{proof}

\begin{lemma}\label{l:positive-NM-1N} Let $(H,\omega)$ be a symplectic Hilbert space defined in \S\ref{s:nullity-Maslov-index}. Let $N\in\Bb(H)$ be an bounded opertor such that $N^*JN=-J$ and $N^2=I_H$. Let $\gamma\in C([a,b],\Sp(H))$ be a positive symplectic path. Then we have $N^*J=-JN$, and the path $\alpha:=N\gamma^{-1}N$ is a positive symplectic path.
\end{lemma}

\begin{proof} Since $N^*JN=-J$ and $N^2=I_H$, we have $N^*J=-JN^{-1}=-JN$. Since $\gamma$ is a symplectic path, $\alpha$ is also a symplectic path. It follows that
\begin{align*}-J\dot\alpha\alpha^{-1}&=-J(N(-\gamma^{-1}\dot\gamma\gamma^{-1})N)(N^{-1}\gamma N^{-1})\\
&=N^*\gamma^*(-J\dot\gamma\gamma^{-1})\gamma N>0.
\end{align*}
\end{proof}

We have the following form of two times iteration formula for the generalized brake symmetry.

\begin{theorem}\label{t:2S-iteration-1} Let $(H,\omega)$ be a symplectic Hilbert space defined in \S\ref{s:nullity-Maslov-index}. Let $N\in\Bb(H)$ be an bounded opertor such that $N^*JN=-J$ and $N^2=I_H$. Let $S\in\Sp(H)$ be such that $(NS)^2=I_H$. Let $\gamma\in C([a,b],\Sp(2n,H))$ be a symplectic path such that $N\gamma(t)^{-1}N\gamma(t)-S$ is Fredholm of index $0$ for each $t\in[a,b]$. Then we have
\begin{align}
\label{e:2S-iteration-1}i_S(N\gamma^{-1}N\gamma)=i_{V^+\times U^+}(\gamma)+i_{V^-\times U^-}(\gamma).
\end{align}
\end{theorem}

\begin{proof} Set $V_0:=\Graph(S)$, $V_1:=V^+\times U^+$, $V_2:=V^-\times U^-$.
By Lemmas \ref{l:index-pm} and \ref{l:2S-iteration}, we have
\[0=\Index(N(\gamma(t))^{-1}N\gamma(t)-S)=\Index(\Graph(\gamma(t)),V_1)+\Index(\Graph(\gamma(t)),V_2)\]
for each $t\in[a,b]$. By \cite[Proposition 1]{BoZh13}, we have $\Index(\Graph(\gamma(t)),V_j)\le 0$ for $j=1,2$. So both of the two index are $0$.
By Lemmas \ref{l:index-pm}, \ref{l:2S-iteration} and \ref{l:positive-NM-1N}, we can apply Theorem \ref{t:nullity-Maslov-index} and our result follows.
\end{proof}

Our next result in the real case was obtained by Y. Long, D. Zhang and the second author \cite[Proposition C]{LoZhZh06} when $S=I_H$, and by X. Hu and S. Sun \cite[Theorem 1.6]{HuSu09} when $S, N\in\Sp(2n,\R)\cap \UUU(2n)$.

\begin{corollary}[Two times iteration formula for the generalized brake symmetry]\label{c:2S-iteration-2} Let $(H,\omega)$ be a complex symplectic Hilbert space of dimension $2n$. Let $N\in\Bb(H)$ be an bounded opertor such that $N^*JN=-J$ and $N^2=I_H$. Let $S\in\Sp(H)$ be such that $(NS)^2=I_H$. Given a path $\gamma\in \Pp_{\tau}(H)$, we have
\begin{align}
\label{e:2S-iteration-2}i_S(\gamma^{(2)})=i_{V^+\times U^+}(\gamma)+i_{V^-\times U^-}(\gamma).
\end{align}
\end{corollary}

\begin{proof} By Lemma \ref{l:Maslov-path-endpoints} and Theorem \ref{t:2S-iteration-1}.
\end{proof}

\Subsection{The iteration formula for the brake symmetry}\label{ss:brake-iteration}

Let $(H,\omega)$ be a symplectic Hilbert space defined in \S\ref{s:nullity-Maslov-index}. Let $N\in\Bb(H)$ be such that $N^2=I_H$ and $N^*JN=-J$. Set $U^{\pm}:=\ker(N\mp I_H)$. Then we have $U^{\pm}\in\Ll(H)$. By replacing the inner product of $H$ with $\lla\cdot,\cdot\rra_{U^-}\oplus\lla\cdot,\cdot\rra_{U^+}$, we can assume that $N^*N=I_H$. Here the form $\omega$ is unchanged. We fix the orthogonal decomposition $H=U^-\oplus U^+$. Then we have $N=\begin{pmatrix}-I_{U^-}& 0\\0& I_{U^+}\end{pmatrix}$, $J=\left(\begin{array}{cc}0&-K^*\\K&0\end{array}\right)$, where $K\in\Bb(U^-,U^+)$ and $K$, $K^*$ are injective maps.

\begin{lemma}\label{l:split-exact-sequence} Let $X$ be a vector space with two linear maps $A,B\in\End(X)$. Assume that one of the following three conditions hold:
\begin{itemize}
\item[(i)] $\ker(AB)\supset\ker A$, and $B:\ker A\to\ker A$ is surjective;
\item[(ii)]  $\ker(AB)\supset\ker A$, $\ker A\cap\ker B=\{0\}$, and $\ker A$ is finite dimensional.
\item[(iii)] $\ker A$ and $\ker B$ are finite dimensional, and $\dim\ker(AB)\ge\dim\ker A+\dim\ker B$.
\end{itemize}
Then there are short exact sequences
\begin{align}\label{e:exact-ker}&0\to\ker B\to\ker(AB)\to\ker A\to 0, \\
\label{e:exact-coker}&0\to X/\ran B\to X/\ran(AB)\to X/\ran A\to 0.
\end{align}
\end{lemma}

\begin{proof} By \cite[Exercises B.11]{Wh78}, it is enough to show that $B:\ker(AB)\to\ker A$ is surjective.

\noi (i) Since $\ker(AB)\supset\ker A$, we have
\[B\ker A\subset B\ker(AB)\subset\ker A.\]
Since $B:\ker A\to\ker A$ is surjective, we have $B\ker(AB)\supset B\ker A=\ker A$. So we obtain $B\ker(AB)=\ker A$. Thus $B:\ker(AB)\to\ker A$ is surjective and our results follow.

\noi (ii) Since $\ker A\cap\ker B=\{0\}$, the map $B:\ker A\to\ker A$ is an injection. Since $\dim\ker A<+\infty$, it is an isomorphism. By (i), our results hold.

\noi (iii) By \cite[Exercises B.11]{Wh78}, we have $\dim\ker(AB)\le\dim\ker A+\dim\ker B$. So we have $\dim\ker(AB)=\dim\ker A+\dim\ker B<+\infty$. By \cite[Exercises B.11]{Wh78}, $B:\ker(AB)\to\ker A$ is an isomorphism. So (\ref{e:exact-ker}) and (\ref{e:exact-coker}) hold.
\end{proof}

\begin{lemma}\label{l:property-NM2=1} Let $M=\begin{pmatrix}A&B\\C&D\end{pmatrix}\in\Sp(H)$ be such that $(NM)^2=I_H$ and $k$ be a positive integer. Then the following hold.
\newline (a) (cf. \cite[Proposition 2.1]{FrOt14}) We have
\begin{align}
\label{e:sym-symplectic1}&KA=D^*K,\;KB=B^*K^*,\;K^*C=C^*K,\\
\label{e:sym-symplectic2}&AB=BD,\;CA=DC,\;A^2-BC=I_{U^-},\;D^2-CB=I_{U^+}.
\end{align}
\newline (b) Let $\a$ be a real number such that $a/\pi\notin\Z$. Then we have
\begin{align}
\label{e:lambda-cos1}\dim \ker(D-(\cos\a)I_{U^+})&=\dim \ker(M-e^{\sqrt{-1}\a}I_H),\\
\label{e:lambda-cos2}\dim (U^-/\ran(D-(\cos\a)I_{U^+}))&=\dim (H/\ran(M-e^{\sqrt{-1}\a}I_H)).
\end{align}
\newline (c) (cf. \cite[Lemma 3.1]{FrOt14}) Denote by $T_k$ and $U_k$ the Chebyshev polynomials of the first kind and the second kind respectively.
Then we have
\begin{equation}\label{e:matrix-power}
M^k=\begin{pmatrix}
T_k(A)& U_{k-1}(A)B\\
CU_{k-1}(A)& T_k(D)
\end{pmatrix}.
\end{equation}
\newline (d) Assume that there is $P\in \Sp(H)$ such that $M=NP^{-1}NP$. Then the pair $(\Graph(M),U^+\times U^+)$ is Fredholm (of index $0$) if  $(\Graph(P),U^+\times U^+)$ and $(\Graph(P),U^+\times U^-)$ are Fredholm (of index $0$). In this case we have
\begin{align}
\label{e:2-iteration-nullity1}\n_{U^+\times U^+}(M)&=\n_{U^+\times U^+}(P)+\n_{U^+\times U^-}(P),\\
\label{e:2-iteration-nullity2}\tilde\n_{U^+\times U^+}(M)&=\tilde\n_{U^+\times U^+}(P)+\tilde\n_{U^+\times U^-}(P).
\end{align}
\newline (e) The pair $(\Graph(M^k),U^+\times U^+)$ is Fredholm (of index $0$) if and only if $(\Graph(M),U^+\times U^+)$ and $M-e^{\sqrt{-1}j\pi/k}I_H$ are Fredholm (of index $0$) for $j=1,\ldots,k-1$. In this case we have
\begin{align}
\label{e:even-iteration-nullity1}\n_{U^+\times U^+}(M^k)&=\n_{U^+\times U^+}(M) + \sum_{j=1}^{k-1}\n_{e^{\sqrt{-1}j\pi/k}}(M),\\
\label{e:even-iteration-nullity2}\tilde\n_{U^+\times U^+}(M^k)&=\tilde\n_{U^+\times U^+}(M) + \sum_{j=1}^{k-1}\tilde\n_{e^{\sqrt{-1}j\pi/k}}(M).
\end{align}
\newline (f) Assume that there is $P\in \Sp(H)$ such that $M=NP^{-1}NP$. Then the pair $(\Graph(PM^k),U^+\times U^+)$ is Fredholm (of index $0$) if and only if $(\Graph(P),U^+\times U^+)$ and $M-e^{2\sqrt{-1}j\pi/(2k+1)}I_H$ are Fredholm (of index $0$) for $j=1,\ldots,k$. In this case we have
\begin{align}
\label{e:odd-iteration-nullity1}\n_{U^+\times U^+}(PM^k)&=\n_{U^+\times U^+}(P) + \sum_{j=1}^k\n_{e^{2\sqrt{-1}j\pi/(2k+1)}}(M),\\
\label{e:odd-iteration-nullity2}\tilde\n_{U^+\times U^+}(PM^k)&=\tilde\n_{U^+\times U^+}(P) + \sum_{j=1}^k\tilde\n_{e^{2\sqrt{-1}j\pi/(2k+1)}}(M).
\end{align}
\end{lemma}

\begin{proof} (a) Since $M$ is symplectic, we have
\[D^*KB-B^*K^*D=0,\;C^*KA-A^*KC=0,\;D^*KA-B^*K^*C=K.\]
Since $M^*JM=J$ and $(NM)^2=I_H$, we have $M^*J=JM^{-1}=JNMN$. Then \eqref{e:sym-symplectic1} holds.
Then we have
\begin{align*}
&K(AB-BD)=D^*KB-B^*K^*D=0,\\
&K^*(CA-DC)=C^*KA-A^*KC=0,\\
&K(A^2-BC)=D^*KA-B^*K^*C=K,\\
&K^*(D^2-CB)=A^*K^*D-C^*KB=K^*.
\end{align*}
By the injectivity of $K$ and $K^*$ we obtain \eqref{e:sym-symplectic2}.
\newline (b) Set $\la:=e^{\sqrt{-1}\a}$. Since $(NM)^2=I_H$, we have $(MN)^2=I_H$ and $\la I_H-MN$ is invertible. Then we have
\begin{align}\label{e:lambda-cos3}
(I_H-\la^{-1}MN)(M-\la I_H)&=M(I_H+N)-(\la I_H+\la ^{-1}N)\\
\nonumber&=\begin{pmatrix}-2 \sqrt{-1} (\sin\a)I_{U^-}&0\\ 0&2(D-(\cos\a) I_{U^+})\end{pmatrix}.
\end{align}
Since $\a/\pi\notin\Z$, \eqref{e:lambda-cos1} and \eqref{e:lambda-cos2} hold.
\newline (c) By the proof of \cite[Proposition 2.1]{FrOt14}, where we replace $A^T$ by $D$.
\newline (d) Let $P$ be of the form
\[P=\begin{pmatrix}A_1&B_1\\C_1&D_1\end{pmatrix},\quad P^{-1}=\begin{pmatrix}A_2&B_2\\C_2&D_2\end{pmatrix}.\]
Then we have $B=2A_2B_1=-2B_2D_1$, and $\ker B_1\cap\ker D_1=\{0\}$.
So we have
\begin{align}
\nonumber\ker B&=MU^+\cap U^+=\{x\in U^+;Px\in NPN(U^+)\}\\
\nonumber &=\{x\in U^+;B_1x=-B_1y, D_1x=D_1y \text{ for some } y\in U^+\}\\
\label{e:ker-B-3} &=\ker B_1\oplus\ker D_1.
\end{align}
By direct calculation, we have
\begin{align*}
  \Graph(M)+U^+\times U^+=&\begin{pmatrix}
                                      I_H&0\\
                                      0&NP^{-1}N\\
                                     \end{pmatrix}
 (\Graph(P)+ U^+\times NPN(U^+))\\
 =&\begin{pmatrix}
                                      I_H&0\\
                                      0&NP^{-1}N\\
                                     \end{pmatrix}
  (U^+\times\{0\} +\\
  &\begin{pmatrix}
  I_{U^-}&0&0\\
  A_1&I_{U^-}&0\\
  C_1&0&I_{U^+}\\
  \end{pmatrix}
 (U^-\times (\ran B_1\oplus\ran D_1)).
 \end{align*}
By Lemma \ref{l:index-pm} we have
\begin{align}\label{e:coker-B-3}
U^-/\ran B\cong H^2/(\Graph(M)+U^+\times U^+)\cong  U^-/\ran B_1\oplus U^+/\ran D_1.
\end{align}
By \eqref{e:ker-B-3} and \eqref{e:coker-B-3}, $B$ is Fredholm if and only if $B_1$ and $D_1$ are Freholm. In this case we have
\begin{align}\label{e:ker-B-2}\dim\ker B&=\dim\ker B_1+\dim\ker D_1,\\
\label{e:coker-B-2}\dim U^-/\ran B&=\dim U^-/\ran B_1+\dim U^+/\ran D_1,\\
\label{e:index-B}\Index B&=\Index B_1+\Index D_1.
\end{align}

By Lemma \ref{l:index-pm} and \eqref{e:co-nullity}, we have $\Index B_1=\Index(\Graph(P),U^+\times U^+)\le 0$, $\Index D_1=\Index(\Graph(P),U^+\times U^-)\le 0$,
and $\Index B=\Index(\Graph(M),U^+\times U^+)\le 0$. So our results follows from \eqref{e:ker-B-2}, \eqref{e:coker-B-2}, and \eqref{e:index-B}.
\newline (e) We have $\ker(U_{k-1}(A)B)\supset \ker B$.
Since $AB=BD$ and $U_{k-1}$ is a polynomial, we have $U_{k-1}(A)B=BU_{k-1}(D)$.
It follows that $\ker(U_{k-1}(A)B)\supset \ker U_{k-1}(D)$.

We recall the identities of Chebyshev polynomials \cite[Lemma 4.1]{FrOt14} for $a\in\CC$,
\begin{equation}\label{e:Chebyshev}
\begin{cases}
T_k(\cos \a)=\cos k\a,\\
U_k(\cos \a)=\frac{\sin ((k+1)\a)}{\sin \a}.
\end{cases}
\end{equation}
It follows that the set of zeros of $U_{k-1}$ is $\{\cos(j\pi/k)| 0< j<k,j\in\Z \}$.
Then we have
\begin{equation}\label{e:UA-dec}
\ker U_{k-1}(D)=\sum_{j=1}^{k-1}\ker(D-\cos(j\pi/k)I_{U^+}).
\end{equation}
By \eqref{e:sym-symplectic2}, we have $D^2-I_n=CB$.
It follows that
\begin{equation}\label{e:ker-B-subset-ker-A}
\ker B\subset \ker (CB)\subset \ker(D-I_{U^+})\oplus \ker (D+I_{U^+}).
\end{equation}
Then we can conclude that \[\ker B\cap \ker (U_{k-1}(D))=\{0\}.\]
By Lemma \ref{l:split-exact-sequence} and (b), $(\Graph(M^k),U^+\times U^+)$ is Fredholm if and only if $(\Graph(M),U^+\times U^+)$ and $M-e^{\sqrt{-1}j\pi/k}I_H$ are Fredholm for $k=1,\ldots,k-1$. In this case \eqref {e:even-iteration-nullity1} and \eqref{e:even-iteration-nullity2} hold. By \eqref{e:co-nullity}, $(\Graph(M^k),U^+\times U^+)$ is Fredholm of index $0$ if and only if $(\Graph(M),U^+\times U^+)$ and $M-e^{\sqrt{-1}j\pi/k}I_H$ are Fredholm of index $0$ for $k=1,\ldots,k-1$.
\newline (f) Let $P$ and $PM^k$ be of the form
\[P=\begin{pmatrix}A_1&B_1\\C_1&D_1\end{pmatrix},\quad PM^k=\begin{pmatrix}A_3&B_3\\C_3&D_3\end{pmatrix}.\]
Since $M=NP^{-1}NP$, we have $PNM=NP$. So $A_1B=B_1(D+I_{U^+})$ holds. Note that $AB=BD$. By (c) we have
\begin{equation}\label{e:B3}B_3=A_1U_{k-1}(A)B+B_1T_k(D)=A_1BU_{k-1}(D)+B_1T_k(D)=B_1R_k(D),\end{equation}
where $R_k(x):=(x+1)U_{k-1}(x)+T_k(x)$ is a $k$-th order polynomial. By \eqref{e:Chebyshev} we have
\begin{equation}R_k(\cos\a)=\frac{\sin ((k+\frac{1}{2})\a)}{\sin \frac{\a}{2}}.\end{equation}
So the the set of zeros of $R_k$ is $\{\cos(2j\pi/(2k+1))| 1\le j\le k,j\in\Z \}$.
Since $M=NP^{-1}NP$ and $NU^+=U^+$, we have $PU^+\cap U^+\subset MU^+\cap PM^kU^+\cap U^+$. So we have $\ker B_1\subset\ker B\cap\ker B_3$. By a proof similar to that of (e), (f) is obtained.
\end{proof}

We have the following form of iteration formula for the brake symmetry.

\begin{theorem}\label{t:brake-iteration} Let $(H,\omega)$ be a symplectic Hilbert space defined in \S\ref{s:nullity-Maslov-index}. Let $N\in\Bb(H)$ be an bounded opertor such that $N^*JN=-J$ and $N^2=I_H$. Let $\gamma\in C([a,b],\Sp(H))$ be a symplectic path
such that $(N\gamma(t))^2=I_H$. Set $U^{\pm}:=\ker(N{\mp}I_H)$. Let $k$ be a positive integer. Then the following hold.
\newline (a) Assume that there is a symplectic path $\gamma_1\in C([a,b],\Sp(H))$ such that $\gamma=N\gamma_1^{-1}N\gamma_1$ and the pair $(\gamma(t)U^+,U^+)$ is Fredholm of index $0$ for each $t\in[a,b]$. Then we have
\begin{align}
\label{e:brake-iteration-21}i_{U^+\times U^+}(\gamma)=i_{U^+\times U^+}(\gamma_1)+i_{U^+\times U^-}(\gamma_1).
\end{align}
\newline (b) Assume that the pair $(\gamma(t)^kU^+,U^+)$ is Fredholm of index $0$ for each $t\in[a,b]$. Then we have
\begin{align}
\label{e:brake-iteration-even1}i_{U^+\times U^+}(\gamma^k)=i_{U^+\times U^+}(\gamma) + \sum_{j=1}^{k-1}i_{e^{\sqrt{-1}j\pi/k}}(\gamma).
\end{align}
\newline (c) Assume that there is a symplectic path $\gamma_1\in C([a,b],\Sp(H))$ such that $\gamma=N\gamma_1^{-1}N\gamma_1$ and the pair $(\gamma_1(t)\gamma(t)^kU^+,U^+)$ is Fredholm of index $0$ for each $t\in[a,b]$. Then we have
\begin{align}
\label{e:brake-iteration-odd1}i_{U^+\times U^+}(\gamma_1\gamma^k)=i_{U^+\times U^+}(\gamma_1) + \sum_{j=1}^k i_{e^{2\sqrt{-1}j\pi/(2k+1)}}(\gamma).
\end{align}
\end{theorem}

\begin{proof} Note that we can assume that $N^*N=I_H$. By Lemma \ref{l:property-NM2=1}, we can apply Theorem \ref{t:nullity-Maslov-index} and our result follows.
\end{proof}

Our next result in the real case was obtained by \cite[Theorem 1.3]{LiZh14}.

\begin{corollary}[The iteration formula for the brake symmetry]\label{c:brake-iteration} Let $(H,\omega)$ be a complex symplectic Hilbert space of dimension $2n$. Let $N\in\Bb(H)$ be an bounded opertor such that $N^*JN=-J$ and $N^2=I_H$. Let $k$ be a positive integer. Given a path $\gamma\in \Pp_{\tau}(H)$, the following hold.
\newline (a) We have
\begin{align}
\label{e:brake-iteration-22}i_{U^+\times U^+}(\gamma^{(2)})=i_{U^+\times U^+}(\gamma)+i_{U^+\times U^-}(\gamma).
\end{align}
\newline (b) Assume that $(N\gamma(\tau))^2=I_H$. Denote by $\tilde\gamma$ the $k$-th $I_H$-iteration of $\gamma$, we have
\begin{align}
\label{e:brake-iteration-even2}i_{U^+\times U^+}(\tilde\gamma)=i_{U^+\times U^+}(\gamma) + \sum_{j=1}^{k-1}i_{e^{\sqrt{-1}j\pi/k}}(\gamma).
\end{align}
\newline (c) We have
\begin{align}
\label{e:brake-iteration-odd2}i_{U^+\times U^+}(\gamma^{(2k+1)})=i_{U^+\times U^+}(\gamma) + \sum_{j=1}^k i_{e^{2\sqrt{-1}j\pi/(2k+1)}}(\gamma^{(2)}).
\end{align}
\end{corollary}

\begin{proof} By Lemma \ref{l:Maslov-path-endpoints} and Theorem \ref{t:brake-iteration}.
\end{proof}

\bigskip

{\bf Acknowledgement}\ \ We would like to thank the referees of this paper for their critical reading and very helpful comments and suggestions.

\end{document}